\documentclass[reqno,a4paper]{amsart}

\usepackage{amssymb}
\usepackage{amsmath}

\makeatletter \@addtoreset{equation}{section} \makeatother

\renewcommand\thefigure{\thesection.\@arabic\c@figure}
\renewcommand\thetable{\thesection.\@arabic\c@table}
\renewcommand\epsilon{\varepsilon}

\newtheorem{theorem}{Theorem}[section]
\newtheorem{lemma}[theorem]{Lemma}
\newtheorem{proposition}[theorem]{Proposition}

\newtheorem{definition}[theorem]{Definition}

\def\eqdef{\stackrel{def}{=}}

\def\emptysquare{{\hbox{\vrule height6pt width0.6pt depth0pt%
\vbox{\hrule height0.6pt width4.8pt depth0pt%
\vglue4.8pt%
\hrule height0.6pt width4.8pt depth0pt}%
\vrule height6pt width0.6pt depth0pt}}}

\def\qed{\unskip\nobreak
\hfil\penalty50\hskip1.75em\null\nobreak\hfil\emptysquare
{\parfillskip=0pt \finalhyphendemerits=0 \par}\medskip}

\begin{document}

\begin{titlepage}

\par\vskip 1cm\vskip 2em

\bigskip

\bigskip

\begin{center}

{\LARGE  {\bf Large deviations for empirical entropies}}\\
{\LARGE  {\bf of $g$-measures}}

\par

\end{center}

\vskip 1cm

\begin{center}

{\LARGE  J.-R. Chazottes$^a$ \ \ D. Gabrielli$^b$}

\par

\end{center}

\bigskip

\bigskip

\begin{center}
$^a$ {\sc Centre de Physique Th{\'e}orique, CNRS-Ecole polytechnique}\\
F-91128 Palaiseau, Cedex, France\\
{\tt jeanrene@cpht.polytechnique.fr}\\
$^b$ {\sc Dipartimento di Matematica Pura e Applicata} \\
Universit{\'a} Dell'Aquila, Via Vetoio Loc. Coppito \\
67100 L'Aquila, Italia\\
{\tt gabriell@univaq.it}
\end{center}

\vskip 1 em \centerline{\bf Abstract}
\smallskip
{\small \noindent } The entropy of an ergodic finite-alphabet
process can be computed from a single typical sample path $x_1^n$
using the entropy of the $k$-block empirical probability and
letting $k$ grow with $n$ roughly like $\log n$. We further assume
that the distribution of the process is a $g$-measure.
We prove large deviation principles
for conditional, non-conditional and relative $k(n)$-block
empirical entropies.

\vfill \noindent {\bf Key words:}\  conditional entropy, relative
entropy, thermodynamic formalism, types.


\bigskip\bigskip

\footnoterule \vskip 1.0em {\small \noindent We acknowledge the
financial support of the Cofin MIUR 2002 prot. 2002027798\_005. We
thank the CIC, Cuernavaca, M{\'e}xico, for its warm hospitality, where
part of this work was done. We also thank A. Galves and D. Guiol
for discussions about the entropy estimation problem.
This project started when both authors were postdoc fellows at IME-USP
(Sa{\~o} Paulo) with a FAPESP grant.
\vskip 1.0em } \noindent

\end{titlepage}

\vfill\eject

%%%%%%%%%%%%%%%%%%%%%% NEW SECTION %%%%%%%%%%%%%%%%%%%%%%%%%%%%%%%%%%%%%%%%%%%%%%%%%%%%%%%%%%%%%%%%%%%%%%%%%%%%%%%%
\section{Introduction}

A problem of interest is the entropy-estimation problem. Given a
sample path $x_{1},x_{2},...,x_{n}$ (where the $x_{i}$'s are drawn
from a finite alphabet $A$) {\em typical} for an unknown ergodic
source, how to estimate its entropy? The simplest idea is to use a
``plug-in'' estimator. First one computes for {\em each} block of
length $k$, the $k$-marginals of the source as the limit, when
$n\to\infty$, of the {\em $k$-block empirical probability} of the
sample $x_1^n$; then one can compute the $k$-block entropy of the
source and let $k\to\infty$ to get the entropy of the source. A
natural question is thus: how is it possible to choose $k=k(n)$ to
do these two steps at the same time?
Ornstein and Weiss \cite{OW} (see also \cite{S}) proved that this is indeed possible for any
ergodic source of positive entropy if $k$ does not grow `too fast'
with $n$, loosely like $\log n$. The proof is based on an
`empirical version' of Shannon-McMillan-Breiman Theorem.

A first result about fluctuations of $k(n)$-block empirical
entropies, refining Ornstein-Weiss' almost-sure result,
was obtained in \cite{3G}. In that paper the
authors consider chains of infinite order which loose memory
exponentially fast. Under additional restrictions on the sequence
$k(n)$ they prove a central limit theorem for the conditional $k(n)$-block
empirical entropy and they prove also that the rescaled $k(n)$-block
empirical entropy cannot have Gaussian fluctuations.

In the present paper, we are interested in large deviations for
$k(n)$-block empirical entropies. To this end we assume
that the distribution of the process generating the sample path $x_1^n$ is
a $g$-measure for the potential $\phi=\log g$ (see below for
definitions and references). Such a process can be viewed as (a
special case of) a chain with complete connections or a
chain of infinite order, see e.g. \cite{FFG,FM1}. Another way,
especially useful for our concern, to characterize and describe a $g$-measure is as a
one-dimensional equilibrium state \cite{FM2,ledrappier}.

In this setting, we prove large deviation principles for
conditional, non-conditional and relative entropies of the
$k(n)$-block empirical probability of the sample path $x_1^n$,
{\em when $k(n)$ grows, roughly speaking, like $\log n$}.
This is done for {\em any} $g$-measure.

When the block length $k$ is fixed, it is easy to obtain a large
deviation principle for $k$-block empirical entropies by
"contraction" of the large deviation principle for the empirical
process \cite{CO}. This is possible because $k$-block entropies
are continuous in the weak topology. To prove the result when
$k(n)$ grows with $n$ we will generalize some classical
combinatorial techniques. We will use the {\em combinatorics of
types} to see "how fast we can let $k$ grow with $n$", and get a
condition close to Ornstein-Weiss' one.

The rate functions we obtain are convex and we will compute
also their Legendre transform which coincide with the corresponding
scaled cumulant generating functions. This will allow to derive some
properties of the rate functions and an explicit representation in
some cases.

Let us notice that the rate function we obtain for conditional and
rescaled non conditional empirical entropy can have a linear part.
This unexpected feature is related to the entropy of
zero-temperature limit of equilibrium states which can be in general
nonzero.

Let us briefly mention that around the problem of
entropy estimation other techniques and ideas have been
developed. The "plug-in" estimator is only one among several
other entropy estimators, see e.g. \cite{CU,CGS,SG,S}.
We point out that we could have worked in the context
of one-dimensional Gibbs measures. An interesting issue
is the case of multi-dimensional Gibbs measures since 
we can no longer use the combinatorics of types.

%%%%%%%%%%%%%%%%%%%%%%%%%%%%%%%%%%%%%%%%%%%%%%%%%%%%%%%%%%%%%%%%%%%%%%%%%%%%%%%%%%%%%%%%%%%%%%%%%%%%%%
The present paper is organized as follows. In the next section we
record preliminary definitions and notions, in particular on
$g$-measures and the various entropies under study. In Section
\ref{main} we present our main results. In Section \ref{comments}
we discuss our results, in particular the form of the rate
functions that we obtain for empirical entropies. Section
\ref{tools} is devoted to the collection of combinatorial tools
needed to understand "how fast $k$ can grow with $n$" later on.
Section \ref{proofs} contains the proof of the main results.

%%%%%%%%%%%%%%%%%%%%%% NEW SECTION %%%%%%%%%%%%%%%%%%%%%%%%%%%%%%%%%%%%%%%%%%%%%%%%%%%%%%%%%%%%%%%%%%%%
\section{Preliminary definitions and notions}\label{setup}

Let $A$ be a finite alphabet. We will denote by
$a_1^\infty\stackrel{def}{=}(a_1,a_2,...)$ the elements of
$A^\mathbb{N}$ and by $a_1^k$ the finite string $(a_1,...,a_k)$.
We will use the notation $x_1^n$ for a ``sample path''
$(x_1,x_2,...,x_n)$, $x_i\in A$. We denote by $T$ the ``shift''
operator defined as $T x_1^\infty = x_2^\infty$. The cylinder set $[a_1^n]$
is the set of infinite strings $b_1^\infty$ drawn from $A^{\mathbb
N}$ such that $b_1^n =a_1^n$.

We call $\mathcal{M}^k$ the set of probability measures $\nu_k$ on
$A^k$ and $\mathcal{M}^k_s$ the set of probability measures
$\nu_k$ on $A^k$ which satisfy the following stationarity condition
\begin{equation}\label{vesuvio}
\sum_{b\in A}\nu_{k}(a_1^{k-1}b)=\sum_{b\in A}\nu_{k}(ba_1^{k-1})
\ \ \ \ \forall a_1^{k-1}\in A^{k-1}\; .
\end{equation}
The subset $\mathcal{M}^k_s$ is convex and $\mathcal{E}^k$ denotes
the set of its extremal elements.

We call $\mathcal{M}$ the set of probability measures $\nu$ on
$A^{\mathbb{N}}$ with the usual sigma-algebra of cylinders. The
subset of shift-invariant (or stationary) measures is denoted by
$\mathcal{M}_s$. The set of ergodic measures (the extremal points
of $\mathcal{M}_s$) is denoted by $\mathcal{E}$.

Given a measure $\nu\in \mathcal{M}_s$ we will write $\nu_k$ for
its $k$-marginals. Of course we have the identity
$\nu_k(a_1^k)=\nu([a_1^k])$ for any $a_1^k\in A^k$ and
consequently $\nu_k\in\mathcal{M}^k_s$.

%%%%%%%%%%%%%%%%%%%%%%%%%%%%%%%%%%%%%%%%%%%%%%%%%%%%%%%%%%%%%%%%%%%%%%%%%%%%%%%%%
\subsection{$g$-measures and equilibrium states}\label{defgmeas}

In this paper we deal with $g$-measures associated to continuous
and regular $g$-functions.
We refer the reader to \cite{keller,ledrappier,walters} for full
details about the following material.

Let $g$ be a continuous function on $A^\mathbb{N}$ satisfying
\begin{equation}\label{defgfunction}
\sum_{b_1^\infty: Tb_1^\infty =a_1^\infty} g(b_1^\infty)
=1\quad\textup{for all}\quad a_1^\infty\in A^{\mathbb N}\,.
\end{equation}
We further assume that $g$ is strictly positive (this implies
$g<1$ by \eqref{defgfunction}). We associate to such a function a
potential, normalized according to \eqref{defgfunction}, by
setting
\begin{equation}\label{defphi}
\phi\stackrel{def}{=} \log g\,.
\end{equation}
Observe that $\phi<0$.
A $g$-measure can be defined as an equilibrium state for the potential $\phi$.
We measure the continuity of $\phi$ by the
sequence of its variations
$(\textup{var}_m(\phi))_{m\in\mathbb{N}}$:
\begin{equation}\label{defvar}
\textup{var}_m(\phi)\stackrel{def}{=}\sup\{|\phi(a_1^\infty)-\phi(b_1^\infty)|:
a_{1}^{m}=b_{1}^{m}\}\,.
\end{equation}
Notice that (uniform) continuity of $\phi$ (with respect to the
canonical distance metrizing product topology) is equivalent to
$\textup{var}_{m}(\phi)\rightarrow 0$ as $m\to\infty$.

It is well-known that if $\textup{var}_m(\phi)$ decreases to $0$
fast enough, then there is a unique $g$-measure which is the
unique equilibrium state for $\phi$.
For instance, if this decreasing is exponential \cite{bowen} or
more generally summable \cite{walters}. 
On another hand, an example of
non-uniqueness was given by Bramson and Kalikow \cite{bramson}. In
that example, $\textup{var}_m(\phi)\geq \frac{C}{\log m}$. Very
recently the authors of \cite{new} showed that square-summability
of variations, ensuring uniqueness \cite{JO}, is tight.
Let us mention a uniqueness criterion based on a
``one-sided'' Dobrushin condition involving oscillations of the
potential instead of variations \cite{FM1}.

From now on, we fix one of the $g$-measures associated to $\phi$
and denote it by $\rho$.
For all $n\geq 1$ and $a_1^\infty\in A^{\mathbb N}$, we have the
following property
\begin{equation}\label{gibbs}
e^{-n \epsilon_n} \leq
\frac{\rho([a_1^n])}{\exp\left(\sum_{j=1}^{n-1}\phi(a_j^\infty)
\right)} \leq e^{n \epsilon_n}
\end{equation}
where $(\epsilon_n)_n$ is a sequence of non-negative real numbers
decreasing to $0$.

For $k\geq 2$, let $\rho^{(k) }$ be the $(k-1)$-step Markov
approximation of $\rho$, that is, the (unique) equilibrium state
of the cylindrical potential
$$
\phi_k(a_1^\infty)=\phi_{k}(a_1^k)\stackrel{def}{=} \log
\frac{\rho([a_1^k])}{\rho([a_2^{k}])}\,\cdot
$$
When $k=1$, $\rho^{(1)}$ is the Bernoulli measure for the
potential $\phi_1(a_1^\infty)=\phi_1(a_1)\stackrel{def}{=} \log
\rho(a_1)$. We can see $\phi_k$ also as a function on $A^k$.

We have the following property
\begin{equation}\label{phi-phi_k}
\parallel \phi -\phi_k \parallel_\infty \leq \textup{var}_k(\phi)\,.
\end{equation}

This implies the statement that for all $a_1^\infty\in
A^{\mathbb{N}}$
$$
\lim_{k\to\infty} \log \frac{\rho([a_1^k])}{\rho([a_2^{k}])}=
\phi(a_1^\infty)
$$
uniformly.

We shall use the variational principle repeatedly. Let
$\psi:A^{{\mathbb N}}\to\mathbb R$ be a continuous function. Then:
\begin{equation}\label{VP}
\sup\{ {\mathbb E}_\eta [\psi] + h(\eta): \eta\in {\mathcal
M}_s\}= P_{top}(\psi)\,.
\end{equation}
Moreover, the supremum is attained if and only if $\eta$ is an
equilibrium state of $\psi$. $P_{top}(\psi)$ is the topological
pressure of $\psi$. It is defined as
\begin{equation}\label{defpressure}
P_{top}(\psi)=\lim_{n\to\infty}\frac{1}{n}\log \sum_{a_1^n}
\exp\big( \sup\big\{\sum_{j=1}^n \psi(b_j^\infty): b\in
[a_1^n]\big\}\big)\,.
\end{equation}

Coming back to a normalized potential $\phi=\log g$, we have
$P_{top}(\phi)=0$. This can be seen, for instance, by plugging
\eqref{gibbs} in \eqref{defpressure}. The variational principle
then tells us that
\begin{equation}\label{entropyccc}
h(\rho)=-\mathbb{E}_{\rho}[\phi]\,.
\end{equation}
In particular, the entropy of a $g$-measure is always strictly
positive.

We shall also consider multiples of the potential $\phi$, that is
potentials of the form $\beta\phi$, $\beta\in\mathbb{R}$. When
$\beta\neq 1$, such potentials have no reason to be normalized as
$\phi$ is, i.e. the corresponding equilibrium states are not
$g$-measures. But this does not matter for us in the sense that we
will only deal with equilibrium states of $\beta\phi$ that we will
indicate with $\rho_{\beta\phi}$.

\bigskip

{\bf Remark}.
A $g$-measure is also named a chain of infinite order or a chain
with complete connections, see e.g. \cite{FFG}, \cite{FM1} for
recent accounts. See also \cite{IO}. In probabilistic terms, a
chain of infinite order, or a chain with complete connections, is
a process characterized by transition probabilities that depend on
the whole past in a continuous manner. A $g$-measure can also be
interpreted as a one-dimensional Gibbs measure if the variations
go to $0$ exponentially fast \cite{FM2}.

%%%%%%%%%%%%%%%%%%%%%%%%%%%%%%%%%%%%%%%%%%%%%%%%%%%%%%%%%%%%%%%%%%%%%%%%%%%%%%%%%%%%%%%%%%%
\subsection{Entropies}\label{entropies}

The $k$-block ($k\geq 1$) Shannon entropy is defined as
\[
H_k(\nu)\stackrel{def}{=}-\sum_{a_{1}^{k}}\nu([a_{1}^{k}])\log\nu([a_{1}^{k}])=
H_k(\nu_k)\stackrel{def}{=}-\sum_{a_{1}^{k}}\nu_{k}(a_{1}^{k})\log\nu_k
(a_{1}^{k})\,.
\]

The conditional $k$-block ($k\geq 2$) entropy is defined as
\[
h_k(\nu)\stackrel{def}{=} -\sum_{a_{1}^{k}}\nu([a_{1}^{k}])\log
\frac{\nu([a_1^k])}{\nu([a_{1}^{k-1}])}
=h_k(\nu_k)\stackrel{def}{=}
-\sum_{a_{1}^{k}}\nu_{k}(a_{1}^{k})\log
\nu_{k}(a_k|a_{1}^{k-1})
\]
where $\nu_k(a_k | a_1^{k-1} )$ is the conditional probability
$\nu_k(a_1^{k})/\sum_b\nu_k( a_1^{k-1}b )$.
We have the relation
\[
h_k(\nu)=H_k(\nu)-H_{k-1}(\nu)\, , k\geq 1
\]
where by convention we set $H_0(\nu)\stackrel{def}{=}0$. Hence
$h_1(\nu)\stackrel{def}{=}H_1(\nu)$. Note that $h_k(\cdot)$ is a
concave function on $\mathcal{M}^k$.

It is well-known that if $\nu$ is a stationary measure, then
\[
\lim_{k\to\infty} h_k(\nu)= \lim_{k\to\infty}
\frac{H_k(\nu)}{k}=h(\nu)
\]
where $h(\nu)$ is the (Shannon-Kolmogorov-Sinai) entropy of $\nu$.

The $k$-block ($k\geq 1$) relative entropy of a stationary measure
$\nu$ with respect to a $g$-measure $\rho$ is defined as
\[
D_k(\nu |\rho) \stackrel{def}{=}
\sum_{a_{1}^{k}}\nu([a_{1}^{k}])\log\frac{\nu([a_{1}^{k}])}{\rho([a_1^k])}
=D_k(\nu_{k}|\rho_{k})\stackrel{def}{=}
\sum_{a_{1}^{k}}\nu_{k}(a_{1}^{k})\log\frac{\nu_k
(a_{1}^{k})}{\rho_k(a_1^k)} \,\cdot
\]
The map $D_k(\cdot | \rho_k)$ is convex on $\mathcal{M}^k$. The
conditional $k$-block ($k\geq 1$) relative entropy is defined as
\[
%\label{relcondkentr}
\Delta_k(\nu|\rho)\stackrel{def}{=}\!\!D_k(\nu|\rho)-D_{k-1}(\nu
|\rho) =\Delta_k(\nu_k |\rho_k)
\stackrel{def}{=}\!\!D_k(\nu_k|\rho_k)-D_{k-1}(\nu_{k-1} |\rho_{k-1})\,.
\]
%We have the relation
%\[
%\label{DDelta}
%\Delta_k(\nu|\rho)=D_k(\nu | \rho) - D_{k-1}(\nu|\rho) \; , k\geq
%1
%\]
Where we set $D_0(\nu|\rho)\stackrel{def}{=}0$. This imposes
$\Delta_1(\nu | \rho)\stackrel{def}{=}D_1(\nu | \rho)$.

The relative entropy $h(\nu|\rho)$ between $\nu\in{\mathcal M}_s$
and a $g$-measure $\rho$ is defined as
\begin{equation}\label{def-relative-entropy}
h(\nu |\rho)\stackrel{def}{=}\lim_{k\to
\infty}\frac{1}{k}D_k(\nu|\rho)= \lim_{k\to
\infty}\Delta_k(\nu|\rho) \quad\textup{and}\quad
h(\nu|\rho)=-\mathbb{E}_{\nu}[\phi]-h(\nu)\,.
\end{equation}
By the variational principle, it is obvious that $h(\nu|\rho)=0$
if, and only if, $\nu$ is an equilibrium state of $\phi$. (See
\cite{CO} for more details.)

\subsection{Empirical measures and entropies}

Given a finite string (a ``sample path'') $x_{1}^{n}$ we define
the empirical measures
\[
\pi_{k}(a_{1}^{k};x_{1}^{n})=\pi_{k,n}(a_{1}^{k})\stackrel{def}{=}
\frac{\sum_{i=1}^{n}{\mathit{1} \!\!\>\!\! I}
(\tilde{x}_{i}^{i+k-1}=a_{1}^{k})}{n} \;\; ,\ k\in \mathbb{N}
\]
where $\tilde{x}_1^\infty\in A^{\mathbb{N}}$ is the periodic, with
period $n$, sample path $(x_{1}^{n}x_{1}^{n}x_{1}^{n} \cdots )$.

It is easy to see that $\pi_{k}(\cdot; x_1^n)\in \mathcal{M}^k_s$.
The family of probability measures $\left(\pi_{k}(\cdot;
x_1^n)\right)_{k\in \mathbb{N}}$ is consistent in the sense that
\[
\sum_{a_j}\pi_{j}(a_1^j;x_1^n)=\pi_{j-1}(a_1^{j-1};x_1^n) \;\; ,\
j\in\mathbb{N}
\]
and are the marginals of the empirical process $\pi(\cdot ;x_1^n)$
defined as
\begin{equation}\label{def-emp-process}
\pi(S ;x_1^n)\stackrel{def}{=}
\frac{1}{n}\sum_{i=1}^{n}\delta_{T^i\tilde{x}_1^\infty}(S)
\end{equation}
where $S$ is any measurable subset of $A^{\mathbb N}$.

\bigskip

We can now define the following plug-in estimators for entropies.

\begin{definition}\label{defentropies}
Let $x_1^n\in A^n$ be a sample path.
The $k$-block empirical entropy is defined as
\[
\hat{H}_k(x_1^n)\stackrel{def}{=}H_k(\pi_{k}(\cdot;x_1^n))\,.
\]
\noindent The conditional $k$-block empirical entropy is defined
as
\[
\hat{h}_k(x_1^n)\stackrel{def}{=}h_k(\pi_{k}(\cdot;x_1^n))\,.
\]
\noindent The relative $k$-block empirical entropy with respect to
a measure $\rho$ is defined as
\[
\hat{D}_k(x_1^n|\rho)\stackrel{def}{=}D_k(\pi_{k}(\cdot;x_1^n)|\rho_k)\,.
\]
\noindent The relative conditional $k$-block empirical entropy
with respect to a measure $\rho$ is defined as
\[
%\label{defCERE}
\hat{\Delta}_k(x_1^n|\rho)\stackrel{def}{=}
\Delta_k(\pi_{k}(\cdot;x_1^n)|\rho_{k})\,.
\]
\end{definition}

%%%%%%%%%%%%%%%%%%%%%%%%%%%%%%%%%%%%%%%%%%%%%%%%%%%%%%%%%%%%%%%%%%%%%%%%%%%%%%%%%%%%%%%%%%%%%%%%%%%%%%%%%%%%%%%%%
\section{Main results}\label{main}

We are now ready to state the main results of this paper.

\begin{theorem}[Large deviation principles for empirical entropies]
\label{primo} Let $x_1^n$ be a sample path distributed according
to a $g$-measure $\rho$. Assume that $\left(k(n)\right)_{n \in
{\mathbb N}}$ diverges and eventually satisfies
\begin{equation}\label{growth}
k(n)\leq \frac{1-\epsilon}{\log|A|}\log n
\end{equation}
for some $0<\epsilon<1$. Then the conditional empirical entropy
$\hat{h}_{k(n)}(x_1^n)$ satisfies the following large deviation
principle:

\noindent For any closed set $C\subset {\mathbb R}$
$$
\limsup_{n\to\infty} \frac{1}{n}\log \rho\left\{ x_1^n :
\hat{h}_{k(n)}(x_1^n) \in C \right\} \leq -\inf\{{\mathbf I}(u):
u\in C\}\,.
$$
For any open set $O\subset {\mathbb R}$
$$
\liminf_{n\to\infty} \frac{1}{n}\log \rho\left\{ x_1^n :
\hat{h}_{k(n)}(x_1^n) \in O \right\} \geq -\inf\{{\mathbf I}(u):
u\in O\}
$$
where the convex rate function ${\mathbf I}$ is defined as
\begin{equation}
{\mathbf I}(u)=\left\{
\begin{array}{cc}\label{I}
\inf\{h(\nu|\rho ):\nu \in \mathcal{M}_s: h(\nu )=u\} & \ \ u\in[0,\log|A|] \\
+\infty & \ \ otherwise\ . \\
\end{array}
\right.
\end{equation}

The same large deviation principle holds if we replace $\hat{h}_{k(n)}(x_1^n)$
by the rescaled empirical entropy $\frac{\hat{H}_{k(n)}(x_1^n)}{k(n)}$.
\end{theorem}

\bigskip

\begin{theorem}[Large deviations for empirical relative entropies]
\label{secondo} Let $x_1^n$ be a sample path distributed according
to a $g$-measure $\rho$. Suppose that $\left(k(n)\right)_{n \in
{\mathbb N}}$ diverges and eventually satisfies $k(n)\leq
\frac{1-\epsilon}{\log|A|}\log n$, for some $0<\epsilon<1$. Then
the empirical relative entropies $\hat{\Delta}_{k(n)}(x_1^n |
\rho)$ and $\frac{1}{k(n)}\hat{D}_{k(n)}(x_1^n|\rho)$ satisfy a
large deviation principle as in Theorem \ref{primo} but with the
rate function
\begin{equation}
{\mathbf J}(u)=\left\{
\begin{array}{cc}
u & \ \ u\in[0,-\inf\{\mathbb{E}_\eta[\phi]:\eta\in\mathcal{E}\}] \\
+\infty & \ \ otherwise\ . \\
\end{array}
\right.
\end{equation}
\end{theorem}

These theorems are proved in Section \ref{proofs}. Their proof relies
in an essential way upon combinatorial properties of types and a
continuity property of entropy which are established in Section \ref{tools}.  

The following proposition deals with the case of fixed block
length. The preceding theorems extend this proposition to
the case when $k(n)$ is allowed to grow with $n$ according to
\eqref{growth}.

\begin{proposition}[Large deviations for fixed block length]\label{kfixed}
Let $x_1^n$ be a sample path distributed according to a
$g$-measure $\rho$. Then, for each $k\geq 1$, the empirical
entropies $\frac{1}{k}\hat{H}_{k}(x_1^n)$, $\hat{h}_{k}(x_1^n)$,
$\frac{1}{k}\hat{D}_{k}(x_1^n)$ and $\hat{\Delta}_{k}(x_1^n)$
satisfy a LDP with normalizing factor $\frac{1}{n}$ and rate
functions respectively given by
$$
{\mathbf I}^H_k(u)=\inf\{h(\nu|\rho ): H_k(\nu)/k=u\}\, ,\quad
{\mathbf I}^h_k(u)=\inf\{h(\nu|\rho ): h_k(\nu )=u\}
$$
$$
{\mathbf I}^D_k(u)=\inf\{h(\nu|\rho ):D(\nu_k|\rho_k)/k=u\}\,
,\quad {\mathbf I}^{\Delta}_k(u)=\inf\{h(\nu|\rho
):\Delta(\nu_k|\rho_k)=u\}
$$
where the infima are taken over $\nu\in \mathcal{M}_s$.
The infimum over an empty set is taken equal to
$+\infty$ following the usual convention.
\end{proposition}

This proposition is a direct consequence of the contraction
principle and suggests that the rate functions we can expect when
we consider $k(n)$ growing with $n$ are ``contracted'' relative
entropies. Note that the rate functions of Proposition
\ref{kfixed} need not be convex.

From the convexity of ${\mathbf I}$ and ${\mathbf J}$ we know that
they are in Legendre duality with the corresponding scaled
cumulant generating function for the different empirical
entropies. In the next two propositions we give the expression of
the scaled cumulant generating function for empirical entropies
and empirical relative entropies.

\begin{proposition}\label{pressk}
Assume that the hypotheses of Theorem \ref{primo} hold. Then the
rate function $\mathbf{I}$ is in Legendre duality with the convex
function $t\mapsto \mathbf{R}(t)$, $t\in\mathbb{R}$, defined as
\begin{equation}
\label{renyibis} {\mathbf R}(t)= \left\{
\begin{array}{l}
(t+1) P_{top}(\phi/(t+1))
\quad\textup{for}\quad t>-1 \\
\sup\{\mathbb{E}_\eta[\phi]:\eta\in\mathcal{E}\}
\quad\textup{for}\quad t\leq -1\,.
\end{array}\right.
\end{equation}
Moreover,
\begin{equation}\label{scgfp}
\lim_{n \to \infty} \frac{1}{n}\log {\mathbb E}_{\rho}\left[
e^{nt\hat{h}_{k(n)}(x_1^n)}\right] =\lim_{n \to \infty}
\frac{1}{n}\log {\mathbb
E}_{\rho}\left[e^{nt\frac{\hat{H}_{k(n)}(x_1^n)}{k(n)}}\right]={\mathbf R}(t)\,.
\end{equation}
\end{proposition}
%%%%%%%%%%%%%%%%%%%%%%%%%%%%%%%%%%%%%%%%%%%%%%%%%%%%%%%%%%%%%%%%%%%%%%%%%%%%%%%%%%%%%%%%%%%%%%%%%

Using \eqref{gibbs} it is easy to check that
\begin{equation}\label{renyi}
{\mathbf R}(t)= (t+1)\ \lim_{n\to\infty} \frac{1}{n}
\log\sum_{a_1^n\in A^n} \rho([a_1^n])^{\frac{1}{t+1}}\quad
\textup{for}\quad t>-1\,.
\end{equation}
This resembles a R{\'e}nyi entropy.

%%%%%%%%%%%%%%%%%%%%%%%%%%%%%%%%%%%%%%%%%%%%%%%%%%%%%%%%%%%%%%%%%%%%%%%%%%%%%%%%%%%%%%%%%%%%%%%%%

%%%%%%%%%%%%%%%%%%%%%%%%%%%%%%%%%%%%%%%%%%%%%%%%%%%%%%%%%%%%%%%%%%%%%%%%%%%%%%%%%%%%%%%%%%%%%%%%%
\begin{proposition}\label{pressrel}
Assume that the hypotheses of Theorem \ref{secondo} hold. Then the
rate function $\mathbf{J}$ is in Legendre duality with the convex
function $t\mapsto \mathbf{P}^{\Delta}(t)$, $t\in\mathbb{R}$,
defined as
$$
{\bf P}^\Delta(t)\stackrel{def}{=} \left\{
\begin{array}{cc}
(1-t)\ \inf\{\mathbb{E}_{\nu}[\phi]:\nu\in \mathcal{E}\} & \ \ t>1 \\
0 & \ \ t\leq 1\,. \\
\end{array}
\right.
$$
Moreover,
\begin{equation}\label{scgfp2}
\lim_{n \to \infty} \frac{1}{n}\log {\mathbb E}_{\rho}\left[
e^{nt\hat{\Delta}_{k(n)}(x_1^n|\rho)}\right] =\lim_{n \to \infty}
\frac{1}{n}\log {\mathbb
  E}_{\rho}\left[e^{nt\frac{\hat{D}_{k(n)}(x_1^n|\rho)}{k(n)}}\right]=
{\mathbf P}^{\Delta}(t)
\end{equation}
\end{proposition}

\bigskip

Let us introduce
\begin{equation}\label{hmin}
h_\infty \stackrel{def}{=} \lim_{\beta\to\infty}
h(\rho_{\beta\phi})\,.
\end{equation}
The existence of this limit will be shown below (Lemma \ref{monotony}).
In general $h_\infty$ can be strictly positive and we stress that it
is equal to $\log |A|$ for the uniform Bernoulli measure.

In the case when $\mathbf{R}$ is a strictly convex, continuously
differentiable function on $]-1,+\infty[$, we can improve the
results of Theorem \ref{primo}. A large class of $g$-measures
satisfies this property, namely those associated to potentials
with square summable variations.

\begin{proposition}[More on large deviations]\label{varphisummable}
In addition to assumptions of Theorem \ref{primo}, assume that the
variations of $\phi$ are square summable.
Then ${\mathbf I}$ is strictly convex on $[h_\infty,\log|A|]$,
with a unique minimum, where it assumes the value $0$, at
$u=h(\rho)$. Moreover it admits the following representation:
\begin{equation}\label{beta}
{\mathbf I}(u)= h(\rho_{\beta_u\phi} |
\rho_\phi)\quad\textup{for}\quad u\in[h_\infty,\log|A|]
\end{equation}
where $\beta_u\geq 0$ is the unique solution of the equation
$h(\rho_{\beta\phi})=u$. On the interval $[0,h_\infty]$ the
function $\mathbf I$ is linear
$$
{\mathbf I}(u)=-u
-\sup\{\mathbb{E}_\eta[\phi]:\eta\in\mathcal{E}\}
$$.
\end{proposition}

%%%%%%%%%%%%%%%%%%%%%%%%%%%%%%%%%%%%%%%%%%%%%%%%%%%%%%%%%%%%%%%%%%%%%%%%%%%%%%%%%%%%%%
\section{Comments on the results}\label{comments}

We make some comments on the above results.

\bigskip

{\bf Zero-temperature limit and non-differentiability of $\mathbf R$
at $-1$.} By using a classical formula for the derivative of the
pressure \cite{keller}, it is straightforward to see that the
right derivative of $t\mapsto \mathbf{R}(t)$ at $-1$, when the
variations of $\phi$ are square summable, is equal to
$$
\lim_{\beta\to+\infty}
\left(P_{top}(\beta\phi)-\beta\mathbb{E}_{\rho_{\beta}}[\phi]\right)
$$
where we recall that $\rho_{\beta\phi}$ is the equilibrium state
of the potential $\beta\phi$.
By the variational principle, we thus get that
$$
\lim_{t\downarrow -1} \frac{d{\mathbf R(t)}}{dt}=
\lim_{\beta\to+\infty} h(\rho_{\beta\phi})=h_{\infty}\,.
$$
This limit is not zero in general, therefore the function
$\mathbf{R}$ is not differentiable at $t=-1$. Notice
that this is related to zero-temperature limit of equilibrium states.

\bigskip

{\bf About the route to large deviations}. Let us emphasize that
we prove our large deviation bounds directly. Another 
way to prove large deviation principles is to first prove the
existence of the corresponding scaled cumulant generating
function, and then to apply G\"artner-Ellis Theorem (see e.g.
\cite{DenH}). To that end one needs to prove, e.g., that the
scaled cumulant generating function is differentiable and strictly
convex. We could do that under the assumption that the potential
of the $g$-measure has square-summable variations. But, as
\eqref{renyibis} shows, the scaled cumulant generating function
is not differentiable at $-1$ in the case $h_{\infty}\neq 0$. Therefore
one cannot apply G\"artner-Ellis theorem. Notice also that the
rate functionals of Proposition \ref{kfixed} can be in general non
convex. This means that even in the case when $k$ is fixed G\"artner-Ellis theorem
may not apply.

We want to stress that with our approach we need not to assume
anything on the rate of convergence to zero of the variations of
the potential.

\bigskip

{\bf On the growth condition \eqref{growth}}. A look at the proof of
Theorem \ref{primo} reveals that we actually have a little bit more
general condition on $k(n)$. In fact we could impose, e.g.,
$$
\frac{(\log n)^2 |A|^{k(n)}}{n}\to 0\quad\textup{as}\;n\to \infty\,.
$$
We feel that condition \eqref{growth} is more appealing and it is
related to the condition which appears in the laws of large numbers
for empirical entropies (see below).

\bigskip

{\bf Flatness of $\mathbf I$.} If $\rho$ is not the unique
equilibrium state of $\phi$, it is easy to see that the rate
function ${\mathbf I}$ can be identically zero in some interval
containing $h(\rho)$. Indeed, the set of equilibrium states of
$\phi$ form a Choquet simplex and the map $\nu\mapsto h(\nu)$ is
convex affine \cite{keller} on the set of shift-invariant
measures. Hence, there is an equilibrium state $\rho_1$ (maybe
equal to $\rho$) such that $h(\rho_1)$ minimizes the entropy among
all equilibrium states of $\phi$. It may be not unique but this
does not matter: we call $h_1$ this minimal entropy. We do the
same for the maximal entropy and call $h_2$ the corresponding
value (maybe equal to $h(\rho)$). Then, it is easy to verify that
${\mathbf I}(u)=0$ for all $u\in[h_1,h_2]$ since ${\mathbf
I}(h_1)={\mathbf I}(h_2)=0$ (by the variational principle) and
${\mathbf I}$ is convex and positive.

\bigskip

{\bf Strong laws of large numbers for empirical entropies.}
If $\rho$ is the unique equilibrium state of $\phi$ (e.g. when $\phi$
has square-summable variations), then $0$ is the minimum of ${\mathbf I}$ and it is attained only at
$u=h(\rho)$ (this is an immediate consequence of the variational
principle). We can use Theorem \ref{primo} and apply Borel-Cantelli Lemma to obtain that
$$
\lim_{n\to +\infty}\frac{1}{k(n)}\hat{H}_{k(n)}(x_1^n)=\lim_{n\to
+\infty}\hat{h}_{k(n)}(x_1^n) =h(\rho)\ \ \ \rho-a.s.
$$
Therefore, we recover in our context the Ornstein-Weiss
almost-sure result cited in the introduction, with a $k(n)$ allowed to grow a little bit
less fast and stronger hypotheses on the source $\rho$. A similar
statement, in probability, can be deduced from the results of
\cite{3G}.
The almost-sure convergence of conditional
empirical entropy in the case of an ergodic measure $\nu$ with
positive entropy can be proved under the condition that $k(n)\leq
\frac{1-\epsilon}{h(\nu)} \log n$ (and $k(n)\to\infty$), for
some $0<\epsilon<1$. If $\epsilon=0$, this almost-sure convergence
fails in general \cite{shields-private}. 

The same argument applied to the statement of Theorem \ref{secondo}
leads to the almost-sure convergence of empirical relative entropies to zero
$$
\lim_{n\to\infty}\hat{\Delta}_{k(n)}(x_1^n|\rho)
=\lim_{n\to\infty}\frac{1}{k(n)}\hat{D}_{k(n)}(x_1^n |\rho)=0 \ \
\ \rho-a.s.
$$
A similar result in probability for $\sqrt{n}\hat{D}_{k(n)}(x_1^n
|\rho)$ appears in \cite{3G} with more assumptions on $k(n)$.

\bigskip

{\bf Connection with central limit asymptotics}.
Theorem \ref{secondo} has its own interest, but it is also
connected with the central limit asymptotics of conditional
empirical entropy \cite{3G} as follows.
The following decomposition holds (see \cite{3G}):
\begin{equation}\label{split}
\hat{h}_{k(n)}(x_1^n) - h(\rho)=
-\frac{1}{n}\sum_{j=0}^{n-1}(\phi(T^j x_1^\infty) - \mathbb{E}_\rho[\phi])
- \hat{\Delta}_{k(n)}(x_1^n|\rho) +\mathcal{C}_n
\end{equation}
where the correction term $\mathcal{C}_n$ is such
that $|\mathcal{C}_n|\leq C\textup{var}_{k(n)}(\phi)$ and $x_1^\infty\in[x_1^n]$.
In words, the conditional empirical entropy is equal to the empirical average of
the potential $-\phi$, plus a term due to the conditional
empirical relative entropy between the empirical measure and the ''true''
measure, and plus a correction.

In \cite{3G}, the authors assume that the variations of $\phi$
decrease exponentially fast. They show, under appropriate
assumptions on the way $k(n)$ is allowed to grow, that
$\sqrt{n}\hat{\Delta}_{k(n)}(x_1^n|\rho)$ goes to zero in
$\rho$-probability, as well as $\sqrt{n}\ \mathcal{C}_n$. Therefore,
they can conclude that the central limit theorem for
$\hat{h}_{k(n)}(x_1^n)-h(\rho)$ is equivalent to the central limit
theorem for $-\frac{1}{n}\sum_{j=0}^{n-1}\phi(T^j x_1^\infty)-
\mathbb{E}_{\rho}[-\phi]$.
In particular, the variance is given by
\begin{equation}\label{variancebis}
\sigma^2=\lim_{n\to\infty} \frac1{n} \mathbb{E}_\rho\big[
\big(\sum_{j=0}^{n-1}\phi(T^j x_1^\infty) -n \mathbb{E}_{\rho}[\phi]\big)^2
\big]\,.
\end{equation}

At large deviation scale it is possible to see that term
$\mathcal{C}_n$ is irrelevant, but not
$\hat{\Delta}_{k(n)}(x_1^n|\rho)$.

In fact large deviations for $\hat{h}_{k(n)}(x_1^n)$ are different
from large deviations for $-\frac{1}{n}\sum_{j=0}^{n-1}\phi(T^j
x_1^\infty)$.
The latter have the same large deviations as
$-\frac{1}{n}\log\rho([x_1^n])$. Indeed, it is easy to check
(using \eqref{gibbs} and \eqref{defpressure}) that for any real
$t$
$$
\Phi(t)\stackrel{def}{=}\!\!\lim_{n\to\infty}
\frac{1}{n}\log\mathbb{E}_{\rho}\left[e^{-t\sum_{j=0}^{n-1}\phi\circ
T^j} \right] \!\!=\!\! \lim_{n\to\infty}\frac{1}{n}\log\sum_{a_1^n\in A^n}
\rho([a_1^n])^{1-t}=\!\!P_{top}((1-t)\phi)\,.
$$
The common rate function for $(-\frac{1}{n}\log\rho([x_1^n]))_n$
and $(-\frac{1}{n}\sum_{j=1}^{n}\phi(T^j x_1^\infty))_n$ is then given
by the Legendre transform of $\Phi$.

In \cite{CP}, it is proved that $\sigma^2=\frac{d^2
\Phi}{dt^2}(0)=\frac{d^2 P_{top}(t\phi)}{dt^2}(0)$. On another
hand, one expects that the second derivative of the scaled
cumulant generating function at $0$ (or, equivalently, the inverse
of the second derivative at $h(\rho)$ of the rate function) equals
the variance (\footnote{Notice that this does not imply a central
limit theorem even under real analyticity, see \cite{bryc}.}).
Though $\mathbf{R}(t)\neq \Phi(t)$ for all $t\neq 0$, a simple
computation shows that $\frac{d^2 \mathbf{R}}{dt^2}(0)=\frac{d^2
P_{top}(t\phi)}{dt^2}(0)=\sigma^2$.

Therefore, we have distinct rate functions (because the
conditional empirical relative entropy ``correction'' contributes
at large deviation scale) but their second derivative at $0$
coincide.

\bigskip

{\bf Remark}.
Using \eqref{split}, the fact that
$\hat{\Delta}_{k(n)}(x_1^n|\rho)\geq 0$, and the fact that
$\mathcal{C}_n$ is irrelevant at large deviation scale, it is easy
to get that
$$
{\mathbf R}(t) \leq \Phi(t)\quad\forall t>0\, ,\quad {\mathbf
R}(t) \geq \Phi(t)\quad\forall t<0\,.
$$

%%%%%%%%%%%%%%%%%%%%%% NEW SECTION %%%%%%%%%%%%%%%%%%%%%%%%%%%%%%%%%%%%%%%%%%%%%%%%%%%%%%%%%%%%%%%%%%%%%%%%%%%%%%%%
\section{Some combinatorial tools}\label{tools}

In this section we collect some definitions and lemmas about types, as
well as a continuity lemma for conditional entropy. These are
essential ingredients for the proofs of our main results which
are in the next section. The proof of the following lemmas
are given in Section \ref{app}

We call $\mathcal{U}^k(A^n)$ the subset of $\mathcal{M}^k_s$ whose
elements can be obtained as empirical measure of sample paths of length n.
Formally we set
\begin{equation}
\mathcal{U}^k(A^n)=\{\nu_k \in \mathcal{M}^k_s : \exists x_1^n\in
A^n\ \ s.t.\ \ \nu_k(\cdot )=\pi_k(\cdot;x_1^n)\}\,.
\end{equation}

The set $A^n$ of sample paths $x_1^n$ can be partitioned into
equivalence classes called types. The equivalence relation
$\sim_k$ is defined as
\begin{equation}
x_1^n\sim_k y_1^n \Leftrightarrow \pi_{k}(\cdot
;x_{1}^{n})=\pi_{k}(\cdot ;y_{1}^{n})\,.
\end{equation}
Let us call $\mathcal{T}^k(A^n)=A^n/\sim_k$ the quotient space.
Elements of $\mathcal{T}^k(A^n)$ are labeled with the
corresponding empirical measure $\pi_{k}(\cdot; x_1^n)$, this
means that there is a bijective correspondence between
$\mathcal{T}^k(A^n)$ and $\mathcal{U}^k(A^n)$. We call
$\tau_{\pi_{k,n}}\in \mathcal{T}^k(A^n)$ the type corresponding to
$\pi_{k}(\cdot; x_1^n)$.

We recall that with ${\mathcal E}^k$ we indicate the extremal
elements of ${\mathcal M}^k_s$.

\begin{lemma}
\label{zerentr} Given a measure $\nu_k \in {\mathcal E}^k$ then
$h_{k}(\nu_k)=0$.
\end{lemma}

\begin{lemma}
\label{addensante} Given a measure $\nu_k\in \mathcal{M}^k_s $
there exists a measure $\mu_k \in \mathcal{U}^k(A^n)$ such that
\begin{equation}
||\mu_k-\nu_k||_{{\textup tv}}=\sum_{a_1^k\in
A^k}|\mu_k(a_1^k)-\nu_k(a_1^k)| \leq \frac{(k+2)|A|^k}{n}
\end{equation}
\end{lemma}

\begin{lemma}\label{shi}
The following inequalities hold
\begin{equation}
\left|\mathcal{T}^k(A^n)\right|=
\left|\mathcal{U}^k(A^n)\right|
\leq (n+1)^{|A|^k} \label{cond1}
\end{equation}
\begin{equation}
|\{x_1^n \in \tau_{\pi_{k,n}}\}|\leq (n-1)e^{nh_{k}(\pi_{k,n})}
\label{bound1}
\end{equation}
\begin{equation}
|\{x_1^n \in \tau_{\pi_{k,n}}\}|\geq (en)^{-2|A|^k}
e^{nh_{k}(\pi_{k,n})} \label{cond3}
\end{equation}
\end{lemma}

\begin{lemma}\label{contlemma}
We have the following continuity property of the conditional
$k$-block entropy:
\begin{equation}
\sup_{\left\{\nu_k,\mu_k:\parallel\nu_k-\mu_k\parallel_{{\textup tv}}\leq
\delta\right\}}\left| h_{k}(\nu_k)-h_{k}(\mu_k)\right|\leq
-2\delta\log\frac{\delta}{|A|^k}
\end{equation}
provided that $\delta\leq e^{-1}$.
\end{lemma}

%%%%%%%%%%%%%%%%%%%%%% NEW SECTION %%%%%%%%%%%%%%%%%%%%%%%%%%%%%%%%%%%%%%%%%%%%%%%%%%%%%%%%%%%%%%%%%%%%%%%%%%%%%%%%
\section{Proofs of main results}\label{proofs}

\subsection{Proof of Theorem \ref{primo}}

Consider a closed set $C\subseteq \mathbb R$. We have
$$
\rho\left\{x_1^n:\ \hat{h}_{k}(x_1^n)\in
C\right\}=\sum_{\left\{x_1^n:\hat{h}_{k}(x_1^n)\in C\right\}}
\rho([x_1^n])\,.
$$

From \eqref{gibbs} and \eqref{phi-phi_k} we get
\begin{equation}
\rho([x_1^n])=e^{n\left\{{\mathbb E}_{\pi_k(\cdot,
x_1^n)}\left[\phi_{k}\right]\right\}} \vartheta_{k,n}(x_1^n)
\label{pro}
\end{equation}
where
\begin{equation}
e^{-n(\epsilon_n + \textup{var}_k(\phi))}\leq
\vartheta_{k,n}(x_1^n) \leq e^{n(\epsilon_n +
\textup{var}_k(\phi))}\,. \label{pro2}
\end{equation}

Hence we have
$$
\sum_{\left\{x_1^n:\hat{h}_{k}(x_1^n)\in C\right\}}
\rho([x_1^n])\leq
$$
$$
e^{n(\epsilon_n + \textup{var}_{k}(\phi))}\times
\sum_{\left\{\pi_{k,n}\in \mathcal{U}^{k}(A^n):h_{k}(\pi_{k,n})\in
C\right\}} \!|\{x_1^n\in \tau_{\pi_{k,n}}\}|\;\;
e^{n\left\{\mathbb{E}_{\pi_{k,n}}\left[\phi_{k}\right]\right\}}
$$
where we have used types defined in Section \ref{tools}. Let us
call
$$
h_k^{-1}(C)\eqdef\left\{\nu_k\in
\mathcal{M}^k_s:h_k(\nu_k)\in C\right\}\quad\textup{and}\quad
h^{-1}(C)\eqdef\left\{\mu\in \mathcal{M}_s:h(\mu)\in C\right\}\,.
$$
Using inequalities \eqref{cond1}-\eqref{bound1} we obtain the
following upper bound
$$
\sum_{\left\{x_1^n:\hat{h}_{k}(x_1^n)\in C\right\}}
\rho([x_1^n])\leq
$$
\begin{equation}
e^{n(\epsilon_n +
\textup{var}_{k}(\phi))}(n+1)^{|A|^k}(n-1)\ \exp\left(
n\left\{\sup_{\nu_{k}\in h_k^{-1}(C)}
(\mathbb{E}_{\nu_{k}}[\phi_{k}]+h_{k}(\nu_{k}))\right\}\right)\,.
\end{equation}
If we consider sequences $(k(n))_{n\in \mathbb N}$ that satisfy the
growth condition \eqref{growth} we obtain
$$
\limsup_{n\to\infty} \frac{1}{n}\log \rho\left\{ x_1^n :
\hat{h}_{k(n)}(x_1^n) \in C \right\}
\leq
\limsup_{k\to \infty}\sup_{\nu_{k} \in h_k^{-1}(C)}
\Big(\mathbb{E}_{\nu_{k}}\left[\phi_{k}\right]+h_{k}(\nu_{k})\Big)\,.
$$
We will prove that for any $\epsilon >0$ there exists an integer $K$ such
that for any $k>K$ and for any $\nu_k\in h_k^{-1}(C)$ there exists
a $\mu\in h^{-1}(C)$ such that
\begin{equation}\label{upperb}
\mathbb{E}_{\nu_{k}}\left[\phi_{k}\right]+h_{k}(\nu_{k})\leq
h(\mu)+\mathbb{E}_{\mu}\left[\phi\right]+\epsilon \,.
\end{equation}
The arbitrariness of $\epsilon$ will imply the first
statement of the theorem.

\noindent To prove formula \eqref{upperb} we have only to take
$\mu$ as the unique $(k-1)$-step Markov extension of $\nu_k$ and
$K$ such that $\textup{var}_{K}(\phi)<\epsilon$.

Let us now prove the lower bound. Consider an open set $O\subseteq
\mathbb R$.
$$
\sum_{\left\{x_1^n:\hat{h}_{k}(x_1^n)\in O\right\}}
\rho([x_1^n])\geq
$$
\begin{equation}\label{down-inf}
e^{-n(\epsilon_n + \textup{var}_{k}(\phi))}\times
\sum_{\left\{\pi_{k,n}\in \mathcal{U}^{k}(A^n):h_{k}(\pi_{k,n})\in
O\right\}} \!|\{x_1^n\in \tau_{\pi_{k,n}}\}|\;\;
e^{n\left\{\mathbb{E}_{\pi_{k,n}}\left[\phi_{k}\right]\right\}}
\end{equation}
Using inequality \eqref{cond3} we obtain
$$
\sum_{\left\{x_1^n:\hat{h}_{k}(x_1^n)\in O\right\}}
\rho([x_1^n])\geq
$$
\begin{eqnarray}
& & e^{-n(\epsilon_n + \textup{var}_{k}(\phi))}(en)^{-2|A|^k}
\sum_{\left\{\pi_{k,n}\in \mathcal{U}^{k}(A^n):h_{k}(\pi_{k,n})\in
O\right\}}e^{n\left\{h_k(\pi_{k,n})+
\mathbb{E}_{\pi_{k,n}}\left[\phi_{k}\right]\right\}}\geq \nonumber \\
& & e^{-n(\epsilon_n + \textup{var}_{k}(\phi))}(en)^{-2|A|^k}
\nonumber
\exp\left(n\left\{\sup_{\left\{\nu_k\in
h_k^{-1}(O)\cap\mathcal{U}^k(A^n)\right\}}\left(h_k(\nu_k)+
\mathbb{E}_{\nu_k}\left[\phi_{k}\right]\right)\right\}\right)\,.
\end{eqnarray}

If we consider sequences $(k(n))_{n\in \mathbb N}$ which satisfy
the growth condition \eqref{growth} we obtain
\begin{eqnarray}
\nonumber
& & \liminf_{n\to\infty} \frac{1}{n}\log \rho\left\{ x_1^n :
\hat{h}_{k(n)}(x_1^n) \in O \right\} \geq \nonumber \\
\nonumber
& &\liminf_{n\to \infty}\sup_{\left\{\nu_{k(n)} \in
h_{k(n)}^{-1}(O)\cap \mathcal{U}^{k(n)}(A^n)\right\}}
\Big(\mathbb{E}_{\nu_{k(n)}}\left[\phi_{k(n)}\right]+h_{k(n)}(\nu_{k(n)})\Big)\,.
\nonumber
\end{eqnarray}
We will prove that for any $\epsilon>0$ and for any $\mu\in
h^{-1}(O)$ there exists a $\pi_{k(n),n}\in h_{k(n)}^{-1}(O)\cap
\mathcal{U}^{k(n)}(A^n)$ such that
$$
\mathbb{E}_{\pi_{k(n),n}}\left[\phi_{k(n)}\right]+h_{k(n)}(\pi_{k(n),n})
\geq h(\mu)+\mathbb{E}_{\mu}\left[\phi\right]-\epsilon\,.
$$
The arbitrariness of $\epsilon$ implies the second
statement of theorem \ref{primo}.

When $n$ is large enough $|h_{k(n)}(\mu_{k(n)})-h(\mu)|$ can become
arbitrarily small and from lemmas \ref{addensante} and
\ref{contlemma}, if $d_n\eqdef (k(n)+2)\frac{|A|^{k(n)}}{n}$, there
exists a measure $\pi_{k(n),n}\in \mathcal{U}^{k(n)}(A^n)$ such
that
$$
\left|h_{k(n)}(\mu_{k(n)})-h_{k(n)}(\pi_{k(n),n})\right|
\leq-2d_n\log\frac{d_n}{|A|^{k(n)}}\,\cdot
$$
For a sequence $(k(n))_{n\in \mathbb N}$ which satisfy the growth
condition \eqref{growth} both $d_n$ and
$-2d_n\log\frac{d_n}{|A|^{k(n)}}$ converge to zero. Since $O$ is an
open set we obtain that if $n$ is large enough there exists a
$\pi_{k(n),n}\in h_{k(n)}^{-1}(O)\cap \mathcal{U}^{k(n)}(A^n))$
and such that $|h_{k(n)}(\pi_{k(n),n})-h(\mu)|$ is arbitrarily
small. It is also easy to show that
$$
\left|\mathbb{E}_{\mu}(\phi)-\mathbb{E}_{\pi_{k(n),n}}(\phi_{k(n)})\right|\leq
\textup{var}_{k(n)}(\phi)\ + d_n\parallel\! \phi\!\parallel_\infty\,.
$$
The statement easily follows.

The proof for the estimator $\frac{\hat{H}_{k(n)}(x_1^n)}{k(n)}$ is
analogous; we will only point out the differences.

For the upper bound we need to prove that for any $\epsilon >0$
there exist a $K$ such that for any $k>K$ and for any $\nu_k\in
\mathcal{M}^k_s$ with $\frac{H_k(\nu_k)}{k}\in C$, there exists
$\mu\in \mathcal{M}_s$ with $h(\mu)\in C$ and such that inequality
\eqref{upperb} holds. This can be done considering
$\mu=\frac{\nu_1^M+\dots +\nu_k^M}{k}$, where $\nu_i^M\in
\mathcal{M}_s$ is the unique $(i-1)$-step Markov extension of
$\nu_i$. Due to the fact that $h$ is affine on $\mathcal{M}_s$, we
have in fact that $h(\mu)=\frac{H_k(\nu_k)}{k}$.

The proof of the lower bound is similar. We omit the details.

The convexity of ${\mathbf I}$ follows from the fact that the maps $h(\cdot
),h(\cdot|\rho):\mathcal{M}_s\to \mathbb R$ are affine.
Given $\nu\in\mathcal{M}_s$ such that $h(\nu)=x$ and
$\mu\in\mathcal{M}_s$ such that $h(\mu)=y$, then for any $c\in
[0,1]$
$$
h(c\nu+(1-c)\mu)=cx+(1-c)y
$$
$$
h(c\nu+(1-c)\mu|\rho)=ch(\nu|\rho)+(1-c)h(\mu|\rho)\,.
$$
This implies that
\begin{equation}
{\mathbf I}(cx+(1-c)y)\leq
h(c\nu+(1-c)\mu|\rho)=ch(\nu|\rho)+(1-c)h(\mu|\rho) \label{cono2}
\end{equation}
If we take the infimum over all $\nu\in\mathcal{M}_s$ such that
$h(\nu)=x$ and $\mu\in\mathcal{M}_s$ such that $h(\mu)=y$ from
\eqref{cono2} one obtains the convexity of ${\mathbf I}$.

Theorem \ref{primo} is proved.

\subsection{Proof of Theorem \ref{secondo}}

The proof of theorem \ref{secondo} is similar to that of Theorem \ref{primo}, so we leave the details
to the reader.

\subsection{Proof of Proposition \ref{kfixed}}

Let us recall the following large deviation principle \cite{CO}.
Let $x_1^n$ be a sample path distributed according to a
$g$-measure $\rho$. Then the empirical process $\pi(\cdot;x_1^n)$
defined at \eqref{def-emp-process} satisfies a large deviation
principle in $(\mathcal{M}_s,d_w)$ with normalizing factor
$\frac{1}{n}$ and rate function
\begin{equation}
I^{\pi}(\nu)=h(\nu|\rho)\,.
\end{equation}
Here $d_w$ is a distance that metrizes weak convergence.

Now we observe that for every fixed $k$ the entropies upon consideration
are continuous in $(\mathcal{M}_s,d_w)$.
Therefore, the contraction principle \cite{DenH} immediately
yields the proposition.

%%%%%%%%%%%%%%%%%%%%%%%%%%%%%%%%%%%%%%%%%%%%%%%%%%%%%%%%%%%%%%%%%%%%%%%%%%%%%%%%%%%%%%%%%%%%%%%%%%%%%%%%
\subsection{Proof of Proposition \ref{pressk}}

We prove that the Legendre transform of $\mathbf I$ is $\mathbf
R$. We know from Theorem \ref{primo} that $\mathbf I$ is a convex
function and this imply the Legendre duality.

We have
%%%%%%%%%%%%%%%%%%%%%%%%%%%%%%%%%%%%%%%%%%%%%%%%%%%%%%%%%%%%%%%%%%%%%%%%%%%%%%%%%%%%%%%%%%%%%
\begin{equation}
\sup_{u\in[0,\log|A|]}\left\{tu\!-\!\!\inf_{\{\nu \in
\mathcal{M}_s: h(\nu )=u\}}h(\nu|\rho )\right\}= \sup_{\nu \in
\mathcal{M}_s}\left\{ {\mathbb E}_{\nu}\left[\phi\right]+
th(\nu)+h(\nu) \right\} \label{romeo}\,.
\end{equation}
If $t> -1$, then we get by applying the variational principle
$$
\eqref{romeo}= (t+1)\sup_{\nu \in \mathcal{M}_s}\left\{{\mathbb
E}_{\nu}\left[\frac{\phi}{t+1}\right]+ h(\nu) \right\}= (t+1)
P_{top}\left(\frac{\phi}{t+1}\right)\,.
$$
It $t< -1$, we get
$$
\eqref{romeo}= (t+1)\inf_{\nu \in \mathcal{M}_s}\left\{{\mathbb
E}_{\nu}\left[\frac{\phi}{t+1}\right]+ h(\nu) \right\}\,.
$$
Observe that $h(\nu)\geq 0$ for all $\nu\in\mathcal{M}_s$.
Moreover, the set of measures with entropy $0$ is dense in
$\mathcal{M}_s$ (wrt weak topology), see e.g. \cite{denker}.
Hence, for $t<-1$,
$\eqref{romeo}=(t+1)\inf\{\mathbb{E}_\eta[\phi/(t+1)] : \eta\in
\mathcal{M}_s\}$. The case $t=-1$ is trivial.

\bigskip

The identification of ${\mathbf R}(t)$ with the scaled cumulant
generating functions (formula \eqref{scgfp}) follows from general
arguments \cite{DenH}.

\bigskip

It is interesting to notice that using the combinatorial properties
of types and the results of Section \ref{tools} it is possible to prove
\eqref{scgfp} directly. We just sketch the proof.

%%%%%%%%%%%%%%%%%%%%%%%%%%%%%%%%%%%%%%%%%%%%%%%%%%%%%%%%%%%%%%%%%%%%%%%%%%%%%%%%%%%%%%%%%%%%%%
%%%%%%%%%%%%%%%%%%%%%%%%%%%%%%%%%%%%%%%%%%%%%%%%%%%%%%%%%%%%%%%%%%%%%%%%%%%%%%%%%%%%%%%%%%%%%%%%

Following arguments already used in the proof of Theorem \ref{primo}
we can obtain
$$
\frac{1}{n}\log \sum_{x_1^n\in A^n} e^{nt\hat{h}_{k(n)}(x_1^n)}\
\rho([x_1^n]) \leq
$$
\begin{equation}\label{up-sup}
\sup_{\nu_{k(n)} \in \mathcal{M}^{k(n)}_s}
\left\{\mathbb{E}_{\nu_{k(n)}}\left[\phi_{k(n)}\right]+(t+1)h_{k(n)}(\nu_{k(n)})\right\}
+ \overline{R}_n
\end{equation}
and
$$
\frac{1}{n}\log \sum_{x_1^n\in A^n} e^{nt\hat{h}_{k(n)}(x_1^n)}\
\rho([x_1^n]) \geq
$$
\begin{equation}\label{down-sup}
\sup_{\nu_{k(n)} \in\mathcal{U}^{k(n)}(A^n)}
\left\{\mathbb{E}_{\nu_{k(n)}}\left[\phi_{k(n)}\right]+(t+1)h_{k(n)}(\nu_{k(n)})\right\}
+ \underline{R}_n
\end{equation}
where $\overline{R}_n$ and $\underline{R}_n$ are correcting terms
converging to zero.

We now compute the supremum in \eqref{up-sup}.

If $t\leq -1$, the function to be maximized is convex and the
supremum is attained at one of the extremal points of
$\mathcal{M}^k_s$, which has entropy zero by virtue of lemma
\ref{zerentr}. Hence the supremum in question equals
\begin{equation}\label{napoli1}
\sup\{\mathbb{E}_{\nu_{k(n)}}[\phi_{k(n)}]: \nu_{k(n)}\in
\mathcal{E}^{k(n)} \}\,.
\end{equation}

If $t>-1$, the supremum in \eqref{up-sup} is equal to
$$
(t+1)\ \sup_{\nu
\in\mathcal{M}_s}\left\{\mathbb{E}_{\nu}\left[\frac{\phi_{k(n)}}{t+1}\right]+h(\nu)
\right\}= (t+1)\ P_{top}\left(\frac{\phi_{k(n)}}{t+1} \right) \,\cdot
$$
To see this, we first notice that if $\nu$ is the $(k(n)-1)$-step
Markov measure having $\nu_{k(n)}$ as $k(n)$-marginals, then
$h_k(\nu_{k(n)})=h(\nu)$. On another hand, the variational
principle tells us that $\mathbb{E}_{\nu}[\phi_{k(n)}/(t+1)] +
h(\nu)$ attains its supremum precisely at a unique $(k(n)-1)$-step
Markov measure because $\phi_{k(n)}$ is a $k(n)$-cylindrical
function. This supremum equals $P_{top}(\phi_{k(n)}/(t+1))$.

It is not difficult to prove now that the limit when $n\to \infty$
of the upper bound coincide with ${\mathbf R}(t)$. Using the
results of Section \ref{tools} it is also possible to prove that
the lower bound has the same limit.

The result for the estimator $\frac{\hat{H}_{k(n)}(x_1^n)}{k(n)}$
can be deduced from the previous result using the fact that $
(\hat{h}_{i}(x_1^n))_i $ is a bounded decreasing sequence and
\begin{equation}\label{average}
\hat{H}_{k}(x_1^n)=\sum_{i=1}^{k} \hat{h}_{i}(x_1^n)\,.
\end{equation}

%%%%%%%%%%%%%%%%%%%%%%%%%%%%%%%%%%%%%%%%%%%%%%%%%%%%%%%%%%%%%%%%%%%%%%%%%%%%%%%%%%%%%%%%%%%%%
%%%%%%%%%%%%%%%%%%%%%%%%%%%%%%%%%%%%%%%%%%%%%%%%%%%%%%%%%%%%%%%%%%%%%%%%%%%%%%%%%%%%%%%%%%%%%%%%%%%%%%%%%%%%

\subsection{Proof of Proposition \ref{pressrel}}

The proof of this proposition is very simple and left to the reader.
It is possible to get \eqref{scgfp2} directly using the combinatorics of
types.

\subsection{Proof of Proposition \ref{varphisummable}}

When the variations of $\phi$ are square summable the map
$\beta\mapsto P_{top}(\beta\phi)$, $\beta\in\mathbb R$, is
continuously differentiable and strictly convex. This can be
deduced from \cite{TV}; The extension of their proofs to the
square summable case is straightforward. This imply that the map
$\mathbf R$ is continuously differentiable and strictly convex in
the interval $(-1,\infty)$. Moreover $\mathbf{R}(0)=0$ and
$\frac{d\mathbf{R}}{dt}(0)=h(\rho)$. This establishes the first
part of the proposition.

\bigskip

We now turn to prove the representation formula \eqref{beta}. First introduce the
following auxiliary function of $\beta\in[0,+\infty)$:
$$
\mathcal{I}(\beta)\stackrel{def}{=}\inf\{ h(\nu|\rho):
\nu\in\mathcal{M}_s, h(\nu)=h(\rho_{\beta\phi})\}\,.
$$
We now claim that $\mathcal{I}(\beta)=h(\rho_{\beta\phi}|\rho)$. The
proof is by contradiction of the variational principle. Assume
that $\eta\neq \rho_{\beta\phi}$ is such that
$$
h(\eta|\rho)\leq h(\rho_{\beta\phi}|\rho)\quad\textup{and}\quad
h(\eta)=h(\rho_{\beta\phi})\,.
$$
This means that (remember \eqref{def-relative-entropy})
$$
\mathbb{E}_\eta[\phi] \geq \mathbb{E}_{\rho_{\beta\phi}}[\phi]\,.
$$
Multiplying this inequality by $\beta>0$ and adding $h(\eta)$ to
the lhs and $h(\rho_{\beta\phi})$ to the rhs (since these two
quantities are indeed equal by hypothesis) yields
$$
\mathbb{E}_\eta[\beta\phi] +h(\eta) \geq
\mathbb{E}_{\rho_{\beta\phi}}[\beta\phi] + h(\rho_{\beta\phi})\,.
$$
But the variational principle tells that the rhs is equal to the
supremum over all shift-invariant measures $\nu$ of
$\mathbb{E}_{\nu}[\beta\phi]+h(\nu)$ and is attained only for
$\nu=\rho_{\beta\phi}$. Therefore $\eta$ must be equal to
$\rho_{\beta\phi}$. In this instance of the variational principle,
we used the fact that if a potential $\phi$ has square summable
variations, then $\beta\phi$ also has square summable variations,
in particular for any $\beta>0$. (\footnote{In case of
non-uniqueness, the claim still holds but $\rho_{\beta\phi}$ is
{\em any} equilibrium state associated to $\beta\phi$ since
relative entropy only depends on $\beta\phi$.})

We now invoke lemma \ref{monotony} hereafter to define a map
$\mathcal{H}: [0,+\infty[\to ]h_\infty,\log|A|]$ defined as
$\mathcal{H}(\beta)=h(\rho_{\beta\phi})$. Since this map is
continuous, strictly decreasing, to each $u\in]h_\infty,\log|A|]$
we can associate a unique $\beta_u$ such that
$h(\rho_{\beta_u})=u$.

The last statement of the proposition follows from the first
comment in Section \ref{comments}.\qed

\bigskip

%%%%%%%%%%%%%%%%%%%%%%%%%%%%%%%%%%%%%%%%%%%%%%%%%%%%%%%%%%%%%%%%%%%%%%%%%%%%%%%%%%%%%%%%%%%%%
We state and prove the lemma used just above.

\begin{lemma}\label{monotony}
Assume that $\phi$ has square summable variations (hence so has
$\beta\phi$ for all $\beta\in\mathbb{R}$) and is not cohomologous
to a constant \textup{(\footnote{I.e. is not the equilibrium
measure for a potential of the form $V-V\circ T +c$, where $V$ is
a measurable function, $c\in\mathbb{R}$. In this case the
equilibrium measure would coincide with the measure of maximal entropy,
the uniform Bernoulli measure.})}. Then the map $\beta\mapsto
h(\rho_{\beta\phi})$ is continuous, strictly decreasing on
$[0,+\infty[$ and $h(\rho_{\beta\phi})\in ]h_\infty,\log |A|]$.
\end{lemma}

\bigskip

{\bf Proof}.
By the variational principle,
$h(\rho_{\beta\phi})=P_{top}(\beta\phi) -
\beta\mathbb{E}_{\rho_{\beta\phi}}[\phi]$. (This shows
continuity.) $\beta\mapsto P_{top}(\beta\phi)$ is strictly
decreasing (since $\phi<0$) and strictly convex (see above). This strict convexity of the pressure can be
translated as follows \cite{keller}
$$
\beta_1 < \beta_2 \quad \Rightarrow\quad \mathbb{E}_{\rho_{\beta_1
\phi}}[\phi] < \mathbb{E}_{\rho_{\beta_2 \phi}}[\phi]\,.
$$
Therefore we get that $\beta\mapsto h(\rho_{\beta\phi})$ is
strictly decreasing when $\beta>0$. It is obvious from the
variational principle that $h(\rho_{\beta\phi})=\log|A|$ when
$\beta=0$. Since $h(\rho_{\beta\phi})$ is bounded from below by
$0$, $h_\infty=\lim_{\beta\to+\infty}h(\rho_{\beta\phi})$ exists.
This ends the proof of the lemma. \qed

%%%%%%%%%%%%%%%%%%%%%%%%%%%%%%%%%%%%%%%%%%%%%%%%%%%%%%%%%%%%%%%%%%%%%%%%%%%%%%%%%%%%%%%%%%%%%%%%%%%%
\section{Proof of some lemmas}\label{app}

This section contains the proof of the lemmas of Section \ref{tools}.

Let us introduce the following graph theoretical representations
that we will use in the proofs of the lemmas.
We call $\mathcal{N}^k_n$ the set of integer-valued maps
$N^k_n:A^k\to \mathbb{N}$ such that
\begin{equation}
\sum_{a_1^k\in A^k}N^k_n(a_1^k)=n
\end{equation}
and
\begin{equation}
\sum_{b\in A}N^k_n(a_1^{k-1}b)=\sum_{b\in A}N^k_n(ba_1^{k-1}) \ \ \
\ \forall a_1^{k-1}\in A^{k-1}\,.
\end{equation}
Let $\mathcal{L}^k_n$ be the subset of $\mathcal{M}^k_s$
whose elements are obtained by normalizing elements in
$\mathcal{N}^k_n$, i.e.
\begin{equation}
\mathcal{L}^k_n=\left\{\nu_k \in \mathcal{M}^k_s: \exists N^k_n\in
\mathcal{N}^k_n \ \ s.t. \ \ \nu_k(\cdot )=\frac{N^k_n(\cdot
)}{n}\right\}\,.
\end{equation}
If $k=1$ then $\mathcal{U}^1(A^n)=\mathcal{L}^1_n$, otherwise a
strict inclusion holds $\mathcal{U}^k(A^n)\subset \mathcal{L}^k_n$
$(n>1)$.

We will call a $k$-order compatible balanced directed multigraph
($k$-multigraph, $k$-M, for short) a directed multigraph with the
following properties: The vertices are labeled with elements of
$A^{k-1}$; For each vertex the number of outgoing arrows is equal to
the number of ingoing arrows; An arrow can go from the vertex
$a_1^{k-1}$ to the vertex $b_1^{k-1}$ if and only if
$a_2^{k-1}=b_{1}^{k-2}$. This arrow inherits the natural label
$a_1^{k-1}b^{k-1}$ (note that several arrows can have the same
label).

Given an element $N^k_n\in \mathcal{N}^k_n$ we represent it with a
$k$-M containing $n$ arrows (\cite{DenH}, section II.2) drawing
$N^k_n(b_1^k)$ directed edges from the vertex associated to
$b_1^{k-1}$ to the one associated to $b_2^k$.

Conversely, given a $k$-M containing $n$ arrows, then it is possible
to associate to it an element of $\mathcal{N}^k_n$ defining
$N_n^k(a_1^{k})$ as the number of arrows going from $a_1^{k-1}$ to
$a_2^{k}$. This gives a bijective correspondence.

To each element $\nu_{k}\in \mathcal{U}^k(A^n)$, we associate the
element $N^k_n=n\nu_{k}\in \mathcal{N}^k_n$.
Then we construct a $k$-M as before, which is connected (note that we are not
considering vertices without ingoing/outgoing arrows).
Given two vertices $a_1^{k-1}$ and $b_1^{k-1}$ which have some 
ingoing/outgoing arrows, there exist $i<j$ with $|i-i|<n$ such that
$\tilde{x}_i^{i+k-2}=a_1^{k-1}$ and $\tilde{x}_j^{j+k-2}=b_1^{k-1}$.
This means that for any $i\leq l<j$ there exists at least one arrow
with label $\tilde{x}_l^{l+k-1}$, i.e., at least one path from the
vertex $a_1^{k-1}$ to the vertex $b_1^{k-1}$.

Conversely given a connected $k$-M we associate to it an element of
$\mathcal{U}^k(A^n)$. A connected $k$-M has at least one Eulerian
circuit (see for example \cite{Bol} section I.3).
One follows the circuit generating a sample path in the following way:
Every time one goes through an arrow with label $a_1^k$, one
concatenates the element $a_k$.
The sample path $x_1^n$ that you obtain in this way is such that
$n\pi(\cdot ;x_1^n)$ has associated the connected $k$-M one started with.

This is a bijective correspondence between $\mathcal{U}^k(A^n)$ and
the subclass of connected $k$-M containing $n$ arrows.
This correspondence says that it is possible, starting from the $k$-M
associated to an element $\pi_{k,n}\in \mathcal{U}^k(A^n)$, to
construct an element $x_1^n\in \tau_{\pi_{k,n}}$ by simply following
an Eulerian circuit.

Some classical combinatorial arguments allow to estimate the number
of Eulerian circuits of a $k$-M and this gives an estimate on the
number of samples $x_1^n \in \tau_{\pi_{k,n}}$ (see \cite{DenH}
section II.2):
\begin{equation}\label{Euler}
\frac{\prod_{a_1^{k-1}}\Big(n\sum_b\pi_{k,n}(a_1^{k-1}b)-1\Big)!}
{\prod_{a_1^{k}}(n\pi_{k,n}(a_1^{k}))!}\leq |\{x_1^n\in
\tau_{\pi_{k,n}}\}| \leq
n\frac{\prod_{a_1^{k-1}}\Big(n\sum_b\pi_{k,n}(a_1^{k-1}b)\Big)!}
{\prod_{a_1^{k}}(n\pi_{k,n}(a_1^{k}))!} 
\end{equation}

We will call a $k$-order weighted compatible balanced directed graph
($k$-weighted graph, $k$-WG, for short) a directed graph with the
following properties: The vertices are labeled with elements of
$A^{k-1}$; To each arrow is associated a nonnegative weight; For
each vertex, the sum of the weights associated to outgoing arrows is
equal to the sum of the weights associated to ingoing arrows; An
arrow can go from the vertex $a_1^{k-1}$ to the vertex $b_1^{k-1}$
if and only if $a_2^{k-1}=b_{1}^{k-2}$; The total sum of the weights
is $1$.

Given a measure $\nu_k\in \mathcal{M}^k_s $ we can represent it by
a $k$-WG and conversely given a $k$-WG we can associate to it an
element of $\mathcal{M}^k_s$.

\subsection{Proof of lemma \ref{zerentr}}

A convex combination of measures corresponds to a $k$-WG with a convex
combination of weights.
Therefore the extremality property in ${\mathcal M}^k_s$ corresponds
to the extremality property in the set of $k$-WG's.
Consider a $k$-WG having nonzero weights only on arrows forming a
single cycle (a loop of successive arrows visiting a vertex no more
than once).
All the nonzero weights are equal to
$\frac{1}{\ell}$, where $\ell$ is the length of the cycle.
Every such a $k$-WG cannot be obtained as a convex combination of other
$k$-WG's. 
Otherwise at least one of them would violate one of the
conditions to be a $k$-WG.
Moreover any $k$-WG can be obtained as a
convex combination of a finite number of $k$-WG's consisting of a
single cycle.
A decomposition can be obtained by iterating a finite
number of times the following procedure.
Take the(an) arrow to which is associated the minimum weight and
consider a cycle containing it.
Substract the minimum weight to all the arrows belonging to the cycle
and add the $k$-WG consisting of the single cycle weighted by
$\frac{m.w.}{\ell}$, where $m.w.=$ minimum weight, to the convex
decomposition.
This gives a complete characterization of ${\mathcal E}^k$.
A direct consequence is that $h_{k}(\nu_k)=0$ for every
$\nu_k\in {\mathcal E}^k$.
This is because for every measure $\nu_k$
with associated a $k$-WG consisting of a single cycle
$\nu_k(a_k|a_1^{k-1})$ can be only zero or one. The lemma is proved.

\subsection{Proof of lemma \ref{addensante}}

Given a measure $\nu_k\in \mathcal{M}^k_s$ it is possible to
construct a measure $\tilde{\mu}_k\in \mathcal{L}^k_n$ such that
$||\tilde{\mu}_k-\nu_k||_{{\textup tv}} \leq \frac{2A^k}{n}$. This
is trivial when $k=1$ and a little bit more tricky when $k>1$
because of the stationarity condition \eqref{vesuvio}. Consider for
any arrow from $a_1^{k-1}$ to $a_2^k$ the following parameter
\begin{equation}
\gamma(a_1^k)=\min\left\{
\left|\nu_k(a_1^k)-\frac{[n\nu_k(a_1^k)]+1}{n}\right|,
\left|\nu_k(a_1^k)-\frac{[n\nu_k(a_1^k)]}{n}\right|\right\}.
\label{gap}
\end{equation}
where $[\cdot ]$ represent the integer part. Take the (an) arrow
with associated the minimum value of $\gamma$. Consider an
elementary cycle containing $a_1^k$ and add or subtract (depending
if the minimum value in \eqref{gap} was obtained with the first or
the second argument) the value $\gamma(a_1^k)$ to all the elements
of the cycle. Fix the values of all the weights whose value is
$\frac{i}{n}$ with $i$ some integer number $0\leq i\leq n$, and
remove them from the $k$-WG. It is easy to see that one can
iterate this procedure up to fix all the values of the weights.
One ends up with some weights which satisfy the stationarity
condition but are not necessarily normalized to one. One concludes
the procedure by adding or subtracting the weight necessary to have
the wanted normalization. One can do this sequentially by using an
elementary unit of weight $\frac{1}{n}$ and adding or subtracting
one unit of weight at the same time in elementary cycles, so that
the stationarity condition is preserved. This is always possible.
The measure $\tilde{\mu}_k$ that is obtained in this way belongs
to $\mathcal{L}^k_n$ and is such that
\[
||\tilde{\mu}_k-\nu_k||_{{\textup tv}}\leq \frac{2|A|^k}{n}\,\cdot
\]

If the $k$-M corresponding to $\tilde{\mu}_k$ is connected then
the proof is finished. If the $k$-M associated to $\tilde{\mu}_k$
is not connected let $m>1$ be the number of connected components
containing respectively $e(1),\cdots ,e(m)$ directed edges with
$\sum_{j=1}^me(j)=n$.
Considering an Eulerian circuit for every component one can associate
a sample path $s(i)$ of length $e(i)$ to the $i$th component
for $i=1,\cdots ,m$.
The measure $\tilde{\mu}_k$ has the following expression
\begin{equation}\label{ping}
\tilde{\mu}_k(\cdot )=\sum_{j=1}^m\frac{e(j)}{n}\pi_k(\cdot ;s(j))\,.
\end{equation}
Let us now consider the sample path $s=s(1)s(2)\cdots s(m)$ of
length $n$ and construct $\mu_k$ as the $k$ empirical measure
$\mu_k(\cdot )=\pi_k(\cdot;s)\in \mathcal{U}^k(A^n)$. Both $\mu_k$
and $\tilde{\mu}_k$ are constructed by sliding windows of width $k$
along cyclicized samples, and computing frequencies in these
windows. Every times the window of size $k$ is overlapped to the
sample $s$ and do not cross points of separation among different
$s(i)$ the $k$-sequence that is matched contributes both in
$\tilde{\mu}_k$ and $\mu_k$.
Using the fact that $m\leq A^{k-1}$ we deduce
\begin{equation}\label{pong}
||\mu_k-\tilde{\mu}_k||_{{\textup tv}}\leq \frac{k|A|^{k-1}}{n}\,\cdot
\end{equation}
Using \eqref{ping}, \eqref{pong} and the triangle inequality yields the statement of the lemma.

\subsection{Proof of lemma \ref{shi} }
The proof of inequalities \eqref{cond1} and \eqref{bound1} is very
simple and elegant and can be found in \cite{S}. More precisely in
section I.6.d in the case of non cyclicized samples and in section
II.1.a in the case of cyclicized samples, which is our case. The
proof of inequality \eqref{cond3} is obtained from estimate
\eqref{Euler} and Stirling formula. Inequality \eqref{bound1} can be
proved in an analogous way.

\subsection{Proof of lemma \ref{contlemma}}

This lemma can be found in \cite{CK} but we give its proof thereafter for the sake of completeness.
Consider $\mu_k$ and $\nu_k$ two measures in $\mathcal{M}^k$ such that
$\parallel\nu_k-\mu_k\parallel_{{\textup tv}}\leq \delta$. Let us set
$\delta_k(a_1^k)=\left|\nu_k(a_1^k)-\mu_k(a_1^k)\right|$.
Obviously
\[
\sum_{a_1^{k-1}}\delta_{k-1}(a_1^{k-1})\leq
\sum_{a_1^k}\delta_k(a_1^k)\leq \delta\,.
\]
Using triangle inequality one obtains
\begin{eqnarray}
\left| h_{k}(\nu_k)-h_{k}(\mu_k)\right| \!\!\!\! &\!\!\leq
\!\! &\!\!\!\!\sum_{a_1^k}\left|\nu_k(a_1^k)\log\nu_k(a_1^k)-\mu_k(a_1^k)\log\mu_k(a_1^k)\right|
\nonumber \\
&+\!\!&\!\!\!\!\sum_{a_1^{k-1}}\left|\nu_{k-1}(a_1^{k-1})\log\nu_{k-1}(a_1^{k-1})
-\mu_{k-1}(a_1^{k-1})\log\mu_{k-1}(a_1^{k-1})\right| \,. \nonumber
\end{eqnarray}
By a simple computation it is possible to obtain the modulus of
continuity of the function $-x\log x$ on the interval $[0,1]$
when $\delta$ is small enough
\[
\sup_{\left\{x,y\in [0,1]:|x-y|\leq
\delta\right\}}\left| x\log x-y\log y\right|=-\delta \log
\delta\,.
\]
Using this result we get
\begin{equation}
\left| h_{k}(\nu_k)-h_{k}(\mu_k)\right|\leq
-\sum_{a_1^k}\delta_k(a_1^k)\log\delta_k(a_1^k)-
\sum_{a_1^{k-1}}\delta_{k-1}(a_1^{k-1})\log\delta_{k-1}(a_1^{k-1})\,.
\label{jensen1}
\end{equation}
We write the right hand side of \eqref{jensen1} as
\begin{equation}
-|A|^k\sum_{a_1^k}\frac{\delta_k(a_1^k)}{|A|^k}\log\delta_k(a_1^k)-
|A|^{k-1}\sum_{a_1^{k-1}}\frac{\delta_{k-1}(a_1^{k-1})}{|A|^{k-1}}
\log\delta_{k-1}(a_1^{k-1})
\end{equation}
and apply Jensen inequality using the fact that $-x\log x$ is a
concave function. When $\delta$ is small enough we finally obtain
\begin{eqnarray}
\left| h_{k}(\nu_k)-h_{k}(\mu_k)\right|&\leq&
-\delta\log\frac{\delta}{|A|^k}-\delta\log\frac{\delta}{|A|^{k-1}}
\nonumber \\
&\leq&-2\delta\log\frac{\delta}{|A|^k}\,\cdot
\end{eqnarray}
The lemma is proved.

\end{document}